\documentclass[a4paper,12pt]{article}
\newtheorem{definition}{{\bf Definition}}[section]
\newtheorem{theorem}[definition]{{\bf Theorem}}
\newtheorem{lemma}[definition]{{\bf Lemma}}

\newtheorem{corollary}[definition]{{\bf Corollary}}

\usepackage{amscd,amsmath,amssymb,graphicx}

\begin{document}

\title{Dynamical entropy of generalized quantum Markov chains
on gauge invariant $C^*$-algebras}
\author{Hiromichi Ohno}
\date{}
\maketitle

\begin{abstract}
We prove that
the mean entropy and the dynamical entropy are equal for
generalized quantum Markov chains on gauge-invariant $C^*$-algebras.

\end{abstract}

\section{Introduction}

The notion of quantum Markov chains was introduced by Accardi
in [\ref{acc}].
As special cases,
the notion of quantum Markov states was defined by Accardi and
Frigerio in [\ref{accardi1}] and
that of $C^*$-finitely correlated states
was discussed by Fannes, Nachtergaele and Werner [\ref{fannes1}].

The notion of generalized quantum Markov chains was
introduced in [\ref{ohno2}, \ref{ohno3}].
Generalized quantum Markov chains 
extend
translation-invariant quantum Markov chains to those on AF algebras.
In [\ref{ohno3}], we considered the case where the AF algebras
are gauge-invariant $C^*$-algebras and
we proved the extendability theorem for
any generalized quantum Markov chain,
that is, the generalized quantum Markov chain is the restriction of
a quantum Markov chain on the UHF algebra.

In [\ref{connes}], Connes, Narnhofer, and Thirring extended the
notion of dynamical entropy of classical dynamical systems to the
case of automorphisms of $C^*$-algebras invariant with respect to
a given state.
We can also find definitions and notations of dynamical entropy 
in [\ref{petz}] for example.

In [\ref{park}],
Park showed that the dynamical entropy of a quantum Markov chain
on the infinite $C^*$-tensor product
 of finite-dimensional $C^*$-algebras is equal to its mean entropy.
In the present paper we discuss the dynamical entropy of 
generalized quantum Markov chains on the gauge-invariant
$C^*$-algebras.

In section 2, we show that the mean entropy of generalized quantum
Markov chains coincides with the mean entropy of extended states on
UHF algebras.
In section 3, we show that the dynamical entropy of a generalized 
quantum Markov chain equal to its mean entropy.


\section{Notation and mean entropy}

Let
${\mathfrak B}_i =M_d = M_d({\mathbb C})$,
the $d \times d$ complex matrix algebra, for $i\in {\mathbb Z}$
and
${\mathfrak B}$ be the infinite $C^*$-tensor product
$ \bigotimes_{i \in {\mathbb Z}} {\mathfrak B}_i$.
For any $i \in {\mathbb Z}$, let
$\beta_i$ be the canonical embedding of $M_d$ to the
$i$th component of 
${\mathfrak B}$.
We denote
${\mathfrak B}_\Lambda =\bigotimes_{n \in \Lambda} {\mathfrak B}_n$ 
for an arbitrary
subset $\Lambda \subset {\mathbb Z}$. The translation $\gamma$ is 
the right shift on ${\mathfrak B}$. We write $\phi_{[1,n]}$ 
for the localization 
$\phi|{\mathfrak B}_{[1,n]}$.
The following definition is from [\ref{fannes1}].


\begin{definition}{\rm
A state $\phi$ on ${\mathfrak B}$ is called a {\it $C^*$-finitely correlated
state} if there exist a finite dimensional $C^*$-algebra ${\mathfrak C}$,
a completely positive map $E:M_d \otimes {\mathfrak C}
\to {\mathfrak C}$ and a state $\rho$ on ${\mathfrak C}$ such that
\[
\rho(E(1_{M_d} \otimes C)) = \rho (C)
\]
for all $C \in {\mathfrak C}$ and
\[
\phi(A_1 \otimes \cdots \otimes A_n)
=\rho(E(A_1 \otimes E(A_2 \otimes \cdots \otimes E(A_n \otimes
1_{\mathfrak C}) \cdots )))
\]
for all $A_1 ,\ldots , A_n \in M_d$.}
\end{definition}

From definition, the $C^*$-finitely correlated state
$\phi$ is translation-invariant.
Therefore, the mean entropy of $\phi$ is defined by
\[
s(\phi)  = \lim_{n\to \infty} {1 \over n} S(D_{[1,n]}),
\]
where $D_{[1,n]}$ is the density matrix of $\phi_{[1,n]}$ and
$S(D_{[1,n]})$ is the von Neumann entropy of $D_{[1,n]}$.

Next, we define generalized quantum Markov chains.
Let $G$ be a unitary subgroup of ${\cal U}(M_d)$.
For any $g\in G$, we can define an automorphism 
$\alpha_g = {\displaystyle\lim_{\longrightarrow}}
{\rm Ad} g^{\otimes n}$ of ${\mathfrak B}$.
We set
\[
{\mathfrak A}_\Lambda = ({\mathfrak B}_\Lambda)^G
=\{ B \in {\mathfrak B}_\Lambda \,\, | \,\, 
\alpha_g ( B )  = B, \,\, g \in G \}
\]
for any finite subset $\Lambda \subset {\mathbb Z}$.
We define
\[
{\mathfrak A} = \overline{\bigcup_{n=1}^\infty {\mathfrak A}_{[-n,n]}}.
\]
Then, ${\mathfrak A}$ is an AF algebra and $\gamma | {\mathfrak A}$ is 
an automorphism of ${\mathfrak A}$.


\begin{definition}
{\rm
Let $\tilde\phi$ be a $C^*$-finitely correlated state
generated by the triple $( M_d , E, \rho)$.
We assume that $E$ satisfies the $G$-covariant condition,
that is,
\begin{eqnarray}
E ( g^{\otimes 2} (A \otimes B) g^{\otimes 2 *} ) =
g E (A\otimes B) g^* \label{eq4.1}
\end{eqnarray}
for all $A,B \in M_d$ and $g \in G$.
Then, we have
\begin{eqnarray}
({\rm id}_{{\mathfrak B}_{[-m,n-1]}} \otimes E) ({\mathfrak A}_{[-m, n+1]})
\subset {\mathfrak A}_{[-m,n]}.\label{eq4.2}
\end{eqnarray}
We call the state $\phi = \tilde\phi |{\mathfrak A}$ a
{\it generalized Markov chain} on ${\mathfrak A}$
generated by $E$ and $\rho$.
}
\end{definition}

From (\ref{eq4.2}),
 we obtain
\[
\phi(A) = \rho \circ E_1 \circ E_2 \circ \cdots \circ E_{n}(A\otimes 1_{M_d})
\]
for any $A \in {\mathfrak A}_{[1,n]}$,
where $E_n = {\rm id}_{{\mathfrak B}_{[1,n-1]}} 
\otimes E| {\mathfrak A}_{[1,n+1]}$.
The above formula justifies the terminology for $\phi$.
Moreover, the extendability theorem in [\ref{ohno3}] says that
any generalized Markov chain (in the sense of [\ref{ohno3}])
on ${\mathfrak A}$ can be written as the restriction of some
$C^*$-finitely correlated state on ${\mathfrak B}$ as above.

In the following, let $\phi$ be the generalized quantum
Markov chain generated by $E$ and $\rho$
and $\tilde\phi$ be the $C^*$-finitely correlated state generated by
the triple $(M_d, E, \rho)$.
This means $\tilde\phi | {\mathfrak A} = \phi$.
Let $\phi_{[1,n]} = \phi | {\mathfrak A}_{[1,n]}$, and
$\tilde{D}_{[1,n]}$ and $D_{[1,n]}$ be the density matrices of 
$\tilde\phi_{[1,n]}$ and $\phi_{[1,n]}$, respectively.

To calculate the mean entropy of $\phi$,
we need next lemma.


\begin{lemma}\label{multi}
Let $\bigoplus_{i=1}^k M_{d_i}$ be a subalgebra of $M_d$ and
$m_i$ be the multiplicity of $M_{d_i}$ in $M_d$.
For any state $\psi$ on $M_d$,
we have
\begin{eqnarray*}
S(D_{\psi | \bigoplus_{i=1}^k M_{d_i}}) 
\le S(D_\psi) + \log k + \log (\max \{m_i \}),
\end{eqnarray*}
where $D$ means the density matrix.
\end{lemma}
\begin{proof}
We consider the extremal decomposition of $\psi$, that is,
\[
\psi = \sum_{j=1}^l \lambda_j \psi_j
\]
for some $\lambda_j \in {\mathbb C}$ and pure states $\psi_j$'s.
Then, we have
\begin{eqnarray*}
S(D_{\psi | \bigoplus_{i=1}^k M_{d_i}}) 
\le \sum_{j=1}^l \eta(\lambda_j) + 
\sum_{j=1}^l \lambda_j S(D_{\psi_j | \bigoplus_{i=1}^k M_{d_i}}),
\end{eqnarray*}
where $\eta (x) = -x \log x$.
Let $F_1$ be a canonical conditional expectation
from $M_d$ to $\bigoplus_{i=1}^k M_{d_i} \otimes M_{m_i}$
and $F_2$ be a conditional expectation from
$\bigoplus_{i=1}^k M_{d_i} \otimes M_{m_i}$ to
$\bigoplus_{i=1}^k M_{d_i}$ using partial traces.
Then, for all $1 \le j \le l$, we have
\[
\dim (F_1(D_{\psi_j})) \le k
\]
and
\[
\dim (F_2 \circ F_1 (D_{\psi_j})) \le k \max \{m_j\}.
\]
This imples
\begin{eqnarray*}
S(D_{\psi | \bigoplus_{i=1}^k M_{d_i}}) 
\le S(D_\psi) + \log k + \log (\max \{m_i \}).
\end{eqnarray*}
\end{proof}

Then, we can get the mean entropy of generalized quantum Markov chains.


\begin{corollary}
We obtain
\[
s(\phi) = s(\tilde\phi).
\]
\end{corollary}
\begin{proof}
It is well known that $s(\tilde\phi) \le s(\phi)$.
Hence, we see the converse.
Let ${\mathfrak A}_{[1,n]} = \bigoplus_{i=1}^{k(n)} M_{d_i}$
and
$m_i(n)$ be the multiplicity of $M_{d_i}$ in ${\mathfrak B}_{[1,n]}$.
From Lemma \ref{multi}, we have
\[
S(D_{[1,n]}) \le S(\tilde{D}_{[1,n]}) + \log k(n) + \log (\max \{m_j(n)\}).
\]
Since $k(n)$ and $m_j(n)$'s have at most 
polynomial growth of $n$, we get
\begin{eqnarray*}
s(\phi) = \lim_{n\to\infty} {1 \over n}S(D_{[1,n]}) 
\le \lim_{n\to\infty} {1 \over n} S(\tilde{D}_{[1,n]})
 = s(\tilde\phi).
\end{eqnarray*}
\end{proof}


\section{Dynamical entropy}
We compute the dynamical entropy following
the work of Park ([\ref{park}]).
First, we remark the well known fact.

\begin{lemma}
Let $F$ be a completely positive map from
$M_d \otimes M_d$ to $M_d$. then there exist
$K_1, \ldots , K_l \in M_d \otimes M_d$ 
such that
\[
F(A) = {\rm Tr}^{(2)}(\sum_{i=1}^l K_i^* A K_i),
\]
where ${\rm Tr}^{(2)}$ is the partial trace over the second component.
\end{lemma}

From this lemma, there exist
$K_1, \ldots , K_l \in M_d \otimes M_d$ 
such that
\[
E(A) = {\rm Tr}^{(2)}(\sum_{i=1}^l K_i^* A K_i).
\]
Let ${\cal D} = {\mathbb C}^l = {\rm span}\{ e_1 , \ldots , e_l\}$
and
$\bar{\mathfrak B}_\Lambda =\bigotimes_{i \in \Lambda}
(M_d \otimes {\cal D})$ for $\Lambda \subset {\mathbb Z}$,
in particular,
$\bar{\mathfrak B} = \bigotimes_{\mathbb Z} 
(M_d \otimes {\cal D})$.
We extend the right shift $\gamma$ to $\bar{\mathfrak B}$ canonically.
Moreover, we extend $\tilde\phi$ to $\bar{\mathfrak B}$ in the
following way.
We define a completely positive map $\bar{E}$ from 
$(M_d \otimes {\cal D}) \otimes M_d$ to $M_d$ by
\[
\bar{E} (A \otimes e_j \otimes B) = {\rm Tr}^{(2)}(K_j^* (A\otimes B) K_j)
\]
for all $A,B\in M_d$ and $1\le j \le l$.
Then, we can define the state $\bar\phi$ on $\bar{\mathfrak B}$ by
\begin{eqnarray*}
&& \bar\phi(A_1\otimes D_1 \otimes A_2 \otimes D_2 \otimes \cdots \otimes
A_{n-1} \otimes D_{n-1} \otimes A_n \otimes D_n) \\
&=& \rho ( \bar{E}(A_1\otimes D_1 \otimes
\cdots \otimes \bar{E}(A_{n-1} \otimes D_{n-1} \otimes
\bar{E}(A_n \otimes D_n \otimes 1_{M_d})) \cdots ))).
\end{eqnarray*}
Furthermore, we put ${K} = \sum_{j=1}^l K_j \otimes e_j \in 
M_d \otimes {\cal D} \otimes M_d$, then we have
\[
\bar{E}= {\rm Tr}^{{\cal D} \otimes M_d}(K^* \, \cdot \, K),
\]
where ${\rm Tr}^{{\cal D} \otimes M_d}$ is a partial trace
from $M_d \otimes {\cal D} \otimes M_d$ to $M_d$.
Let $\alpha_1$ and $\alpha_2$ be canonical embeddings of
${\mathfrak A}$ to ${\mathfrak B}$ and of ${\mathfrak B}$ to 
$\bar{\mathfrak B}$, respectively.
From a simple calculation, we have
\[
\tilde\phi = \bar\phi \circ \alpha_2
\]
and
\[
\phi = \bar\phi \circ \alpha_2 \circ \alpha_1.
\]


Now, we consider the dynamical entropy of the generalized quantum Markov 
chain $\phi$.
Since it is well known that
$h_\phi(\gamma) \le s(\phi)$, we show the converse.
From the property of dynamical entropy, we have
\begin{eqnarray}
&&h_\phi(\gamma) = 
\lim_{n\to\infty} h_{\phi,\gamma}({\mathfrak A}_{[1,n]}) \nonumber \\
&\ge& \lim_{n\to\infty} \lim_{k\to\infty}{1\over kn}
H_\phi ({\mathfrak A}_{[1,n]}, {\mathfrak A}_{[n+1, 2n]}, \ldots ,
{\mathfrak A}_{[(k-1)n+1,nk]}) \nonumber \\
&\ge& 
\lim_{n\to\infty} \lim_{k\to\infty}{1\over kn}
H_{\tilde\phi} (\alpha_1({\mathfrak A}_{[1,n]}),
 \alpha_1({\mathfrak A}_{[n+1, 2n]}), \ldots ,
\alpha_1({\mathfrak A}_{[(k-1)n+1,nk]})) \nonumber \\
&\ge&
\lim_{n\to\infty} \lim_{k\to\infty}{1\over kn}
H_{\bar\phi} (\alpha_2\circ\alpha_1({\mathfrak A}_{[1,n]}),
 \ldots ,
\alpha_2\circ\alpha_1({\mathfrak A}_{[(k-1)n+1,nk]})) \label{first}.
\end{eqnarray}


To consider the last term, we make a decomposition of $\bar\phi$.
Let 
\[
{\cal I} = \{ 1, \ldots ,d \}^n \times \{1 , \ldots , l\}^n 
\]
and
\[
{\cal J} = \{ J = (J_1,J_2, \ldots J_k) \,|\,
J_i \in {\cal I} \}.
\]
For any $J_i = (a_1, \ldots , a_n) \times (b_1,\ldots,b_n) \in {\cal I}$,
 we set
\[
P_{J_i} = (P_{a_1} \otimes e_{b_1}) \otimes \cdots \otimes (P_{a_n} \otimes
e_{b_n}) \in \bar{\mathfrak B}_{[1,n]},
\]
where $P_1, \ldots , P_d$ are mutually orthogonal 
minimal projections in $M_d$.
Moreover, we put
\[
P_J = P_{J_1} \otimes P_{J_2} \otimes \cdots \otimes P_{J_k}
\in \bar{\mathfrak B}_{[1,nk]}
\]
for any $J \in {\cal J}$.
Let $W \in M_d$ be the density matrix of $\rho$ and
\[
K_{[0,nk+1]} = \delta_{nk}(K)\cdot 
\delta_{nk-1}(K) \cdots 
\delta_{0}(K),
\]
where $\delta_{i}$ is the canonical embedding of 
$M_d \otimes {\cal D} \otimes M_d$ to its copy in $\bar{\mathfrak B}$
located at the $i$th ${\cal D}$.
Now, we define a linear functional
 $\bar\phi_J$ on $\bar{\mathfrak B}_{[1,nk]}$ for any
$J \in {\cal J}$ by
\begin{eqnarray}\label{phiJ}
\bar\phi_J = {\rm Tr} (K_{[0,nk+1]} P_J \beta_0(W) 
K_{[0,nk +1]}^* \,\, \cdot \,\,).
\end{eqnarray}
Since $P_J$ and $\beta_0(W)$ commute, $\bar\phi_J$ is positive.
Moreover, we have
\[
\bar\phi_{[1,nk]} = \sum_{J \in {\cal J}} \bar\phi_J .
\]
Similarly, we can define the positive linear functional 
$\bar\phi_I$ on $\bar{\mathfrak B}_{[1,n]}$ for
any $I \in {\cal I}$


By using this decomposition, we have
\begin{eqnarray}
&& H_{\bar\phi} (\alpha_2\circ\alpha_1({\mathfrak A}_{[1,n]}),
 \ldots ,
\alpha_2\circ\alpha_1({\mathfrak A}_{[(k-1)n+1,nk]})) \nonumber \\
&\ge&
\sum_{J \in {\cal J}} \eta(\bar\phi_J(1)) -k \sum_{I \in {\cal I}}
\eta(\bar\phi_I(1)) \nonumber \\
&& + k S(\phi_{[1,n]}) - k \sum_{I \in {\cal I}}
\bar\phi_I(1) S(\widehat{\bar\phi}_I | {\mathfrak A}_{[1,n]}),
\label{second}
\end{eqnarray}
where $\widehat\cdot$ means the normalization of 
positive linear functionals.


Let ${\cal M} = {\rm span}\{ P_1, \ldots , P_d\} 
\otimes {\cal D}$.
We can define the classical Markov chain $\mu$ on 
$\bigotimes_{\mathbb Z} {\cal M}$ generated by
$\bar{E}|{\cal M} \otimes {\cal M}$ and $\rho |{\cal M}$.
Then, we have
\[
\mu (P_J) = \bar\phi_J(1),
\]
for all $J \in {\cal J}$. Therefore,
\[
\lim_{n\to\infty} \lim_{k\to\infty}
{1\over kn} \sum_{J\in {\cal J}} \eta(\bar\phi_J(1))
\]
converges to the mean entropy of $\mu$.
We obtain 
\begin{eqnarray} 
\lim_{n\to\infty} \lim_{k\to\infty}
{1\over kn}(\sum_{J \in {\cal J}} \eta(\bar\phi_J(1)) -k \sum_{I \in {\cal I}}
\eta(\bar\phi_I(1)) ) =0. \label{third} 
\end{eqnarray}


On the other hand, by Lemma \ref{multi},
\[
S(\widehat{\bar\phi}_I | {\mathfrak A}_{[1,n]})
\le  S(\widehat{\bar\phi}_I | {\mathfrak B}_{[1,n]}) + \log P(n)
\]
for some polynomial $P$.
Hence, we compute $S(\widehat{\bar\phi}_I | {\mathfrak B}_{[1,n]})$.
We consider $\widehat{\bar\phi}_I$ as a state on 
$\bar{\mathfrak B}_{[0,n]} \otimes M_d$ in the way of (\ref{phiJ}).
Since the density matrix of $\widehat{\bar\phi}_I$ is 
$K_{[0,n+1]} P_I \beta_0(W) 
K_{[0,n +1]}^*$ and $P_I$ is a minimal projection
in $\bar{\mathfrak B}_{[1,n]}$, we can decompose 
\[
\widehat{\bar\phi}_I = \sum_{j=1}^{ld^2} \lambda_j \omega_j
\]
for some $\lambda_j \in {\mathbb C}$ and pure states $\omega_j$'s
on $\bar{\mathfrak B}_{[0,n]} \otimes M_d$.
Since $\omega_j$ is a pure state on ${\mathfrak B}_{[0,n+1]} \otimes 
\bigotimes_{i=0}^n {\cal D}$,
we can write
\[
\omega_j = \omega_{j,1} \otimes \omega_{j,2},
\]
where $\omega_{j,1}$ and $\omega_{j,2}$ are pure states
on ${\mathfrak B}_{[0,n+1]}$ and $\bigotimes_{i=0}^n {\cal D}$,
respectively.
Therefore, we obtain
\begin{eqnarray*}
&&S(\widehat{\bar\phi}_I | {\mathfrak B}_{[1,n]}) \\
&\le& S(\widehat{\bar\phi}_I | {\mathfrak B}_{[0,n+1]})
+S(\widehat{\bar\phi}_I | {\mathfrak B}_{\{0,n+1\}}) \\
&\le& \sum_{j=1}^{ld^2} \lambda_j S(\omega_j | {\mathfrak B}_{[0,n+1]})
+\sum_{j=1}^{ld^2} \eta(\lambda_j) + 2 \log d \\
&=& \sum_{j=1}^{ld^2} \lambda_j S(\omega_{j,1} )
+\sum_{j=1}^{ld^2} \eta(\lambda_j) + 2 \log d \\
&\le& 4 \log d + \log l.
\end{eqnarray*}
Now, we have
\begin{eqnarray}
&& \lim_{n\to\infty} \lim_{k \to \infty} {1 \over kn}
k \sum_{I \in {\cal I}}
\bar\phi_I(1) S(\widehat{\bar\phi}_I | {\mathfrak A}_{[1,n]}) \nonumber \\
&\le& \lim_{n\to\infty}{1 \over n} 
(\log P(n) + 4\log d + \log l )
=0. \label{fourth}
\end{eqnarray}
Combining (\ref{first}), (\ref{second}),
(\ref{third}) and (\ref{fourth}), we obtain
\[
h_\phi(\gamma) \ge s(\phi).
\]
Consequently, we have the next theorem.
\begin{theorem}
Let $\phi$ be a generalized quantum Markov chain on ${\mathfrak A}$.
Then, we obtain
\[
h_\phi(\gamma) = s(\phi).
\]
\end{theorem}


\end{document}